\documentclass[12pt,reqno]{amsart}
\usepackage{amsfonts}
\usepackage{amsmath,amsthm,amsfonts,amssymb}
\newtheorem{lemma}{Lemma}[section]
\newtheorem{theorem}{Theorem}[section]
\newtheorem{corollary}{Corollary}[section]

\def\bl{\begin{lemma}}
\def\bt{\begin{theorem}}
\def\el{\end{lemma}}
\def\et{\end{theorem}}
\def\bp{\begin{proof}}
\def\ep{\end{proof}}
\def\bc{\begin{corollary}}
\def\ec{\end{corollary}}

\def\mb{\mathbb}

\def\-{\setminus}

\def\vp{\varphi}

\def\lt{\left}
\def\rt{\right}
\def\+{\bigcup}
\def\.{\bigcap}

\def\ll{\langle}
\def\rl{\rangle}

\title[A right inverse of Cauchy-Riemann operator $\overline{\partial}^k+a$]
{A right inverse of Cauchy-Riemann operator $\overline{\partial}^k+a$ in weighted Hilbert space $L^2(\mb{C},e^{-|z|^2})$}

\author[Shaoyu Dai and Yifei Pan]{Shaoyu Dai$^1$ and Yifei Pan$^2$}

\address{1 School of Mathematics, Southeast University, Nanjing, 210096, China;}
\address{
Department of Mathematics, Jinling Institute of Technology, Nanjing, 211169, China.}
\address{\it E-mail address: dymdsy@163.com}

\address{2 Department of Mathematical Sciences, Purdue University Fort Wayne, Fort Wayne, 46805-1499, USA.}
\address{\it E-mail address: pan@pfw.edu}

\begin{document}

\begin{abstract} Using H\"{o}rmander $L^2$ method for Cauchy-Riemann
equations from complex analysis, we study a simple differential operator $\overline{\partial}^k+a$ of any order (densely defined and closed) in weighted Hilbert space $L^2(\mb{C},e^{-|z|^2})$ and prove
the existence of a right inverse that is bounded.
\end{abstract}

\maketitle

\section{Introduction}

In this paper, using H\"{o}rmander $L^2$ method \cite{1} for Cauchy-Riemann
equations from complex analysis, we study the right inverse of the differential operator $\overline{\partial}^k+a$, which is densely defined and closed, in a Hilbert space by proving the following result on the existence of (entire) weak solutions of the equation $\overline{\partial}^ku+au=f$ in the weighted Hilbert space $L^2(\mb{C},e^{-|z|^2})$. Here and throughout, $a$ is a complex constant, $k$ a positive integer, $\overline{\partial}^k:=\frac{\partial^k}
{\partial\overline{z}^k}$, $k^{th}$-order Cauchy-Riemann operator, where $\frac{\partial}{\partial\overline{z}}=\frac{1}{2}
(\frac{\partial}{\partial_x}+i\frac{\partial}{\partial_y})$, and $d\sigma:=\frac{1}{2i}d\overline{z}\wedge dz$, the volume form.

\bt\label{thzy}
For each $f\in L^2(\mb{C},e^{-|z|^2})$, there exists a weak solution $u\in L^2(\mb{C},e^{-|z|^2})$ solving the equation
$$\overline{\partial}^ku+au=f$$ in $\mb{C}$
with the norm estimate $$\int_\mb{C} |u|^2e^{-|z|^2}d\sigma\leq
\frac{1}{k!}\int_\mb{C} |f|^2e^{-|z|^2}d\sigma.$$
\et

The novelty of Theorem \ref{thzy} is that the differential operator $\overline{\partial}^k+a$ has a bounded (linear) right inverse
\begin{align*}
T_k: L^2(\mb{C},e^{-|z|^2})&\longrightarrow L^2(\mb{C},e^{-|z|^2}), \\
(\overline{\partial}^k+a)T_k&=I
\end{align*}
with the norm estimate $\|T_k\|\leq\frac{1}{\sqrt{k!}}$.
In particular, the differential operator $\overline{\partial}^k$ has a bounded right inverse $T: L^2(\mb{C},e^{-|z|^2})\longrightarrow L^2(\mb{C},e^{-|z|^2})$, which, to the best of our knowledge, appears to be new. We also note the fact that the constant $a$ dose not appear in the norm estimate, and it is this fact that we shall use later.

For the first order $\overline{\partial}:=\overline{\partial}^1$, the Cauchy-Riemann operator, we have the following slight extension of the simplest case of H\"{o}rmander's theorem in the complex plane (\cite{3} and \cite{4}) ($a=0$; see \cite{5} for a related result). Note that $\Delta=4\partial\overline{\partial}$.

\bt\label{th1}
Let $\vp$ be a smooth and nonnegative function on $\mb{C}$ with
$\Delta\vp>0$. For each $f\in L^2(\mb{C},e^{-\vp})$ such that $\frac{f}{\sqrt{\Delta\vp}}\in L^2(\mb{C},e^{-\vp})$, there exists a weak solution $u\in L^2(\mb{C},e^{-\vp})$ solving the equation
$$\overline{\partial}u+au=f$$
with the norm estimate
$$\int_\mb{C} |u|^2e^{-\vp}d\sigma\leq
4\int_\mb{C} \frac{|f|^2}{\Delta\vp}e^{-\vp}d\sigma.
$$
\et

The organization of the paper is as follows. In Section 2, we will prove several key lemmas based on functional analysis in terms of H\"{o}rmander $L^2$ method, while the proof of Theorem \ref{thzy} and \ref{th1} will be given in Section 3. In Section 4, we will give some further remarks.

\section{Several lemmas}

Here, we consider the weighted Hilbert space
$$L^2(\mb{C},e^{-\vp})
=\{f\mid f\in L^2_{loc}(\mb{C}); \int_\mb{C}|f|^2 e^{-\vp}d\sigma<+\infty\},$$
where $\vp$ is a nonnegative function on $\mb{C}$.
We denote
the weighted inner product for $f,g\in L^2(\mb{C},e^{-\vp})$ by
$\langle f,g\rangle_\vp=\int_\mb{C}\overline{f}g e^{-\vp}d\sigma$,
and the weighted norm of $f\in L^2(\mb{C},e^{-\vp})$ by
$\|f\|_\vp=\sqrt{\langle f,f\rangle_\vp}.$
Let $C_0^\infty(\mb{C})$ denote the set of all smooth functions $\phi: \mb{C}\rightarrow\mb{C}$ with compact support. For $u,f\in L^2_{loc}(\mb{C})$, we say that $f$ is the $k^{th}$ weak $\overline{\partial}$ partial derivative of $u$, written $\overline{\partial}^ku=f$, provided $\int_\mb{C}u\overline{\partial}^k\phi d\sigma=(-1)^k\int_\mb{C}f\phi d\sigma$ for all test functions $\phi\in C_0^\infty(\mb{C})$; we say that $f$ is the $k^{th}$ weak $\partial$ partial derivative of $u$, written $\partial^ku=f$, provided $\int_\mb{C}u\partial^k\phi d\sigma=(-1)^k\int_\mb{C}f\phi d\sigma$ for all test functions $\phi\in C_0^\infty(\mb{C})$.

In the following, let $\vp$ be a smooth and nonnegative function on $\mb{C}$. For $\forall\phi\in C_0^\infty(\mb{C})$, we first define the following formal adjoint of $\overline{\partial}^k$ with respect to the weighted inner product in $L^2\lt(\mb{C},e^{-\vp}\rt)$. Let $u\in L^2_{loc}(\mb{C})$. We integrate as follows by the definition of the weak partial derivative.
\begin{align*}
\lt\ll \phi,\overline{\partial}^ku\rt\rl_\vp&=\int_\mb{C}\overline{\phi}\
\lt(\overline{\partial}^ku\rt) e^{-\vp}d\sigma \\
&=(-1)^k\int_\mb{C}\lt(\overline{\partial}^k\lt(\overline{\phi} e^{-\vp}\rt)\rt)ud\sigma\\
&=(-1)^k\int_\mb{C}e^{\vp}\lt(\overline{\partial}^k\lt(\overline{\phi} e^{-\vp}\rt)\rt)ue^{-\vp}d\sigma\\
&=(-1)^k\int_\mb{C}\overline{e^{\vp}\partial^k\lt(\phi e^{-\vp}\rt)}ue^{-\vp}d\sigma\\
&=\lt\ll(-1)^ke^{\vp}\partial^k\lt(\phi e^{-\vp}\rt),u\rt\rl_\vp\\
&=:\lt\ll \overline{\partial}_\vp^{k*}\phi, u\rt\rl_\vp,
\end{align*}
where $\overline{\partial}_\vp^{k*}\phi =(-1)^ke^{\vp}\partial^k\lt(\phi e^{-\vp}\rt)$ is so called the formal adjoint of $\overline{\partial}^k$ with domain in $C_0^\infty(\mb{C})$. Let $\lt(\overline{\partial}^k+a\rt)_\vp^{*}$ be the formal adjoint of $\overline{\partial}^k+a$ with domain in $C_0^\infty(\mb{C})$. Note that $I_\vp^{*}=I$, where $I$ is  the identity operator. Then $\lt(\overline{\partial}^k+a\rt)_\vp^{*}=\overline{\partial}_\vp^{k*}+a$.

Now we give several lemmas for a general weight based on functional analysis, which are the core elements of H\"{o}rmander $L^2$ method.

\bl\label{lemmaifif}
For each $f\in L^2(\mb{C},e^{-\vp})$, there exists a entire weak solution $u\in L^2(\mb{C},e^{-\vp})$ solving the equation
$$\overline{\partial}^ku+au=f$$ in $\mb{C}$
with the norm estimate
$$\|u\|^2_\vp\leq c$$
if and only if
$$|\langle f,\phi\rangle_\vp|^2\leq c\lt\|\lt(\overline{\partial}^k+a\rt)^*_\vp\phi\rt\|^2_\vp, \ \ \forall\phi\in C_0^\infty(\mb{C}),$$
where $c$ is a constant.
\el

\bp
Let $\overline{\partial}^k+a=H$. Then $\lt(\overline{\partial}^k+a\rt)^*_\vp=H^*_\vp$.

(Necessity) For $\forall\phi\in C_0^\infty(\mb{C})$, from the definition of $H^*_\vp$ and Cauchy-Schwarz inequality, we have
$$
|\langle f,\phi\rangle_\vp|^2=|\langle Hu,\phi\rangle_\vp|^2=\lt|\lt\langle u,H^*_\vp\phi\rt\rangle_\vp\rt|^2\leq\|u\|^2_\vp
\lt\|H^*_\vp\phi\rt\|^2_\vp\leq c\lt\|H^*_\vp\phi\rt\|^2_\vp
=c\lt\|\lt(\overline{\partial}^k+a\rt)^*_\vp\phi\rt\|^2_\vp.
$$

(Sufficiency) Consider the subspace
$$E=\lt\{H^*_\vp\phi\mid\phi\in C_0^\infty(\mb{C})\rt\}\subset L^2(\mb{C},e^{-\vp}).$$
Define a linear functional $L_f: E\rightarrow\mb{C}$ by
$$L_f\lt(H^*_\vp\phi\rt)=\langle f,\phi\rangle_\vp=\int_\mb{C}\overline{f}\phi e^{-\vp}d\sigma.$$
Since
$$\lt|L_f\lt(H^*_\vp\phi\rt)\rt|=\lt|\langle f,\phi\rangle_\vp\rt|
\leq\sqrt{c}\lt\|H^*_\vp\phi\rt\|_\vp,$$
then $L_f$ is a bounded functional on $E$. Let $\overline{E}$ be the closure of $E$ with respect to the norm $\|\cdot\|_\vp$ of $L^2(\mb{C},e^{-\vp})$. Note that $\overline{E}$ is a Hilbert subspace of
$L^2(\mb{C},e^{-\vp})$. So by Hahn-Banach's extension theorem, $L_f$ can be extended to a linear functional $\widetilde{L}_f$ on $\overline{E}$
such that
\begin{equation}\label{29}
\lt|\widetilde{L}_f(g)\rt|\leq\sqrt{c}\lt\|g\rt\|_\vp, \ \ \forall g\in \overline{E}.
\end{equation}
Using the Riesz representation theorem for $\widetilde{L}_f$, there exists a unique $u_0\in \overline{E}$ such that
\begin{equation}\label{30}
\widetilde{L}_f(g)=\langle u_0,g\rangle_\vp, \ \ \forall g\in \overline{E}.
\end{equation}

Now we prove $\overline{\partial}^ku_0+au_0=f$. For $\forall\phi\in C_0^\infty(\mb{C})$, apply $g=H^*_\vp\phi$ in (\ref{30}). Then
$$\widetilde{L}_f\lt(H^*_\vp\phi\rt)=\lt\langle u_0,H^*_\vp\phi\rt\rangle_\vp=\lt\langle Hu_0,\phi\rt\rangle_\vp.$$
Note that
$$\widetilde{L}_f\lt(H^*_\vp\phi\rt)=L_f\lt(H^*_\vp\phi\rt)=\langle f,\phi\rangle_\vp.$$
Therefore,
$$\lt\langle Hu_0,\phi\rt\rangle_\vp=\langle f,\phi\rangle_\vp, \ \ \forall\phi\in C_0^\infty(\mb{C}),$$
i.e.,
$$\int_\mb{C}\overline{Hu_0}\phi e^{-\vp}d\sigma=\int_\mb{C}\overline{f}\phi e^{-\vp}d\sigma, \ \ \forall\phi\in C_0^\infty(\mb{C}).$$
Thus, $Hu_0=f$, i.e., $\overline{\partial}^ku_0+au_0=f$.

Next we give a bound for the norm of $u_0$. Let $g=u_0$ in (\ref{29}) and (\ref{30}). Then we have
$$\|u_0\|^2_\vp=\lt|\langle u_0,u_0\rangle_\vp\rt|=\lt|\widetilde{L}_f(u_0)\rt|
\leq\sqrt{c}\lt\|u_0\rt\|_\vp.$$
Therefore, $\|u_0\|¡ª_\vp^2\leq c$.

Note that $u_0\in\overline{E}$ and $\overline{E}\subset L^2(\mb{C},e^{-\vp})$. Then $u_0\in L^2(\mb{C},e^{-\vp})$. Let $u=u_0$. So there exists $u\in L^2(\mb{C},e^{-\vp})$ such that
$\overline{\partial}^ku+au=f$ with $\|u\|^2_\vp\leq c$. The proof is complete.
\ep

\bl\label{lemmaH}
\begin{align*}
\lt\|\lt(\overline{\partial}^k+a\rt)^*_\vp\phi\rt\|^2_\vp
=\lt\|\lt(\overline{\partial}^k+a\rt)\phi\rt\|^2_\vp+\lt\langle \phi,\overline{\partial}^k\lt(\overline{\partial}_\vp^{k*}\phi\rt)-
\overline{\partial}^{k*}_\vp
\lt(\overline{\partial}^k\phi\rt)\rt\rangle_\vp, \ \ \forall\phi\in C_0^\infty(\mb{C}).
\end{align*}
\el

\bp
Let $\overline{\partial}^k+a=H$. Then $\lt(\overline{\partial}^k+a\rt)^*_\vp=H^*_\vp$.
For $\forall\phi\in C_0^\infty(\mb{C})$,
\begin{align}
\lt\|H^*_\vp\phi\rt\|^2_\vp&=\lt\langle H^*_\vp\phi,H^*_\vp\phi\rt\rangle_\vp\nonumber\\
&=\lt\langle \phi,HH^*_\vp\phi\rt\rangle_\vp\nonumber\\
&=\lt\langle \phi,H^*_\vp H\phi\rt\rangle_\vp+\lt\langle \phi,HH^*_\vp\phi-H^*_\vp H\phi\rt\rangle_\vp\nonumber\\
&=\lt\langle H\phi,H\phi\rt\rangle_\vp+\lt\langle \phi,HH^*_\vp\phi-H^*_\vp H\phi\rt\rangle_\vp\nonumber\\
&=\lt\|H\phi\rt\|^2_\vp+\lt\langle \phi,HH^*_\vp\phi-H^*_\vp H\phi\rt\rangle_\vp.\label{31}
\end{align}
Note that
\begin{align*}
HH^*_\vp\phi=\lt(\overline{\partial}^k+a\rt)\lt(\overline{\partial}^k+a\rt)^*_\vp\phi
=\overline{\partial}^k\lt(\overline{\partial}_\vp^{k*}\phi\rt)+a\overline{\partial}^k\phi+a\overline{\partial}_\vp^{k*}\phi+a^2\phi
\end{align*}
and
\begin{align*}
H^*_\vp H\phi&=\lt(\overline{\partial}^k+a\rt)^*_\vp\lt(\overline{\partial}^k+a\rt)\phi
=\overline{\partial}^{k*}_\vp\lt(\overline{\partial}^k\phi\rt)+a\overline{\partial}_\vp^{k*}\phi+a\overline{\partial}^k\phi+a^2\phi.
\end{align*}
Then
\begin{equation}\label{32}
HH^*_\vp\phi-H^*_\vp H\phi=\overline{\partial}^k\lt(\overline{\partial}_\vp^{k*}\phi\rt)-\overline{\partial}^{k*}_\vp\lt(\overline{\partial}^k\phi\rt).
\end{equation}
So by (\ref{31}) and (\ref{32}), we have
\begin{align*}
\lt\|H^*_\vp\phi\rt\|^2_\vp=\lt\|H\phi\rt\|^2_\vp+\lt\langle \phi,\overline{\partial}^k\lt(\overline{\partial}_\vp^{k*}\phi\rt)-\overline{\partial}^{k*}_\vp\lt(\overline{\partial}^k\phi\rt)\rt\rangle_\vp.
\end{align*}
This lemma is proved.
\ep

\bl\label{lemmak} For $\forall\phi\in C_0^\infty(\mb{C})$,
\begin{align}\label{4}
\overline{\partial}^k\lt(\overline{\partial}_\vp^{k*}\phi\rt)
-\overline{\partial}^{k*}_\vp\lt(\overline{\partial}^k\phi\rt)
&=(-1)^k\sum^k_{i=1}
\sum^k_{j=1}\lt(\begin{array}{c}
  k  \\
  i
\end{array}\rt)\lt(\begin{array}{c}
  k  \\
  j
\end{array}\rt)\partial^{k-i}\overline{\partial}^{k-j}\phi \overline{\partial}^jP_i,
\end{align}
where
\begin{equation}\label{27}
P_i=\sum\frac{i!}{m_1!m_2!\cdots m_i!}
\prod^i_{\gamma=1}\lt(\frac{-\partial^\gamma\vp}{\gamma!}\rt)^{m_\gamma},
\end{equation}
and the sum is over all $i$-tuples of nonnegative integers $(m_1,m_2,\cdots m_i)$ satisfying the constraint $1m_1+2m_2+\cdots+ im_i=i$.
\el

\bp
\begin{align}
\overline{\partial}_\vp^{k*}\phi&=(-1)^ke^{\vp}\partial^k\lt(\phi e^{-\vp}\rt)\nonumber\\
&=(-1)^ke^{\vp}\sum^k_{i=0}\lt(\begin{array}{c}
  k  \\
  i
\end{array}\rt)\partial^{k-i}\phi \partial^ie^{-\vp}\nonumber\\
&=(-1)^k\sum^k_{i=0}\lt(\begin{array}{c}
  k  \\
  i
\end{array}\rt)\lt(\partial^{k-i}\phi\rt)\lt(e^{\vp}
\partial^ie^{-\vp}\rt)\label{1}.
\end{align}
Then from (\ref{1}), we have
\begin{align*}
\overline{\partial}^k\lt(\overline{\partial}_\vp^{k*}\phi\rt)
&=(-1)^k\sum^k_{i=0}\lt(\begin{array}{c}
  k  \\
  i
\end{array}\rt)\overline{\partial}^k
\lt(\lt(\partial^{k-i}\phi\rt)\lt(e^{\vp}\partial^ie^{-\vp}\rt)\rt)\\
&=(-1)^k\sum^k_{i=0}\lt(\begin{array}{c}
  k  \\
  i
\end{array}\rt)
\left(\sum^k_{j=0}\lt(\begin{array}{c}
  k  \\
  j
\end{array}\rt)\overline{\partial}^{k-j}\partial^{k-i}\phi \overline{\partial}^j\lt(e^{\vp}\partial^ie^{-\vp}\rt)\right)\\
&=(-1)^k\sum^k_{i=0}
\sum^k_{j=0}\lt(\begin{array}{c}
  k  \\
  i
\end{array}\rt)\lt(\begin{array}{c}
  k  \\
  j
\end{array}\rt)\partial^{k-i}\overline{\partial}^{k-j}\phi \overline{\partial}^j\lt(e^{\vp}\partial^ie^{-\vp}\rt)
\end{align*}
and
\begin{align*}
\overline{\partial}^{k*}_\vp\lt(\overline{\partial}^k\phi\rt)
=(-1)^k\sum^k_{i=0}\lt(\begin{array}{c}
  k  \\
  i
\end{array}\rt)\lt(\partial^{k-i}\overline{\partial}^k\phi\rt)
\lt(e^{\vp}\partial^ie^{-\vp}\rt).
\end{align*}
Therefore,
\begin{align}\label{2}
\overline{\partial}^k\lt(\overline{\partial}_\vp^{k*}\phi\rt)
-\overline{\partial}^{k*}_\vp\lt(\overline{\partial}^k\phi\rt)
=(-1)^k\sum^k_{i=1}
\sum^k_{j=1}\lt(\begin{array}{c}
  k  \\
  i
\end{array}\rt)\lt(\begin{array}{c}
  k  \\
  j
\end{array}\rt)\partial^{k-i}\overline{\partial}^{k-j}\phi \overline{\partial}^j\lt(e^{\vp}\partial^ie^{-\vp}\rt).
\end{align}
Let $h(g)=e^g, g=-\vp$. By Fa\`{a} di Bruno's formula \cite{2},
\begin{align}
\partial^ie^{-\vp}&=\partial^i(h(g))\nonumber\\
&=\sum\frac{i!}{m_1!m_2!\cdots m_i!}h^{(m_1+\cdots+m_i)}(g)
\prod^i_{\gamma=1}\lt(\frac{\partial^\gamma g}{\gamma!}\rt)
^{m_\gamma}\nonumber\\
&=\lt(\sum\frac{i!}{m_1!m_2!\cdots m_i!}
\prod^i_{\gamma=1}\lt(\frac{-\partial^\gamma\vp}{\gamma!}\rt)^{m_\gamma}\rt)
e^{-\vp}\nonumber\\
&:=P_ie^{-\vp},\label{3}
\end{align}
where the sum is over all $i$-tuples of nonnegative integers $(m_1,m_2,\cdots m_i)$ satisfying the constraint $1m_1+2m_2+\cdots+ im_i=i$.
So (\ref{4}) is proved by (\ref{2}) and (\ref{3}).
\ep

However, unlike in the previous lemmas, here we have to be confined with a special weight $\vp=|z|^2$ in order to deal with the high order differential operator and we note that the space $L^2(\mb{C},e^{-|z|^2})$ is a well-known space, sometimes, called Fock space.
\bl\label{lemmak2}
Let $\vp=|z|^2$. Then
\begin{align}\label{28}
\lt\ll\phi,\overline{\partial}^k\lt(\overline{\partial}_\vp^{k*}\phi\rt)
-\overline{\partial}^{k*}_\vp\lt(\overline{\partial}^k\phi\rt)
\rt\rl_\vp=\sum^{k-1}_{j=0}\frac{(k!)^2}{(j!)^2(k-j)!}\lt\| \overline{\partial}^{j}\phi\rt\|^2_\vp, \ \ \forall\phi\in C_0^\infty(\mb{C}).
\end{align}
\el

\bp
By $\vp=|z|^2$, we have
\begin{align*}
\partial^\gamma\vp=\left\{
\begin{array}{ccc}
\overline{z},       &      & {\gamma=1,}\\
0,     &      & {\gamma\geq2.}
\end{array} \right.
\end{align*}
Then from (\ref{27}),
\begin{align}\label{9}
P_i=\lt(-\partial\vp\rt)^i=\lt(-\overline{z}\rt)^i=(-1)^i\lt(\overline{z}\rt)^i.
\end{align}
Note that
\begin{align*}
\overline{\partial}^j\overline{z}^i=\left\{
\begin{array}{ccc}
0,       &      & {i<j,}\\
\frac{i!}{(i-j)!}\overline{z}^{i-j},     &      & {i\geq j}.
\end{array} \right.
\end{align*}
Let $s=i-j$.
So for $\forall\phi\in C_0^\infty(\mb{C})$, by (\ref{4}) and (\ref{9}) we have
\begin{align}
\overline{\partial}^k\lt(\overline{\partial}_\vp^{k*}\phi\rt)
-\overline{\partial}^{k*}_\vp\lt(\overline{\partial}^k\phi\rt)
&=(-1)^k\sum^k_{j=1}
\sum^k_{i=j}\lt(\begin{array}{c}
  k  \\
  i
\end{array}\rt)\lt(\begin{array}{c}
  k  \\
  j
\end{array}\rt)\lt(\partial^{k-i}\overline{\partial}^{k-j}\phi\rt)
\frac{(-1)^i i!}{(i-j)!}\overline{z}^{i-j}\nonumber\\
&=(-1)^k\sum^k_{j=1}
\sum^{k-j}_{s=0}\lt(\begin{array}{c}
  k  \\
  j+s
\end{array}\rt)\lt(\begin{array}{c}
  k  \\
  j
\end{array}\rt)\lt(\partial^{k-j-s}\overline{\partial}^{k-j}\phi\rt)
\frac{(-1)^{j+s}(j+s)!}{s!}\overline{z}^s\nonumber\\
&=(-1)^k\sum^{k}_{j=1}A_{k-j}\sum^{k-j}_{s=0}\lt(\begin{array}{c}
  k-j  \\
  s
\end{array}\rt)\lt(\partial^{k-j-s}\overline{\partial}^{k-j}\phi\rt)
(-1)^s\overline{z}^s\nonumber\\
&=(-1)^k\sum^{k}_{j=1}A_{k-j}\sum^{k-j}_{s=0}\lt(\begin{array}{c}
  k-j  \\
  s
\end{array}\rt)\lt(\partial^{k-j-s}\overline{\partial}^{k-j}\phi\rt)
P_s,\label{6}
\end{align}
where $$A_{k-j}=\frac{(-1)^j(k!)^2}{((k-j)!)^2j!}.$$
Then by (\ref{6}) and (\ref{3}), we have
\begin{align*}
\lt(\overline{\partial}^k\lt(\overline{\partial}_\vp^{k*}\phi\rt)
-\overline{\partial}^{k*}_\vp\lt(\overline{\partial}^k\phi\rt)
\rt)e^{-\vp}
&=(-1)^k\sum^{k}_{j=1}A_{k-j}\sum^{k-j}_{s=0}\lt(\begin{array}{c}
  k-j  \\
  s
\end{array}\rt)\lt(\partial^{k-j-s}\overline{\partial}^{k-j}\phi\rt)
P_se^{-\vp}\\
&=(-1)^k\sum^{k}_{j=1}A_{k-j}\sum^{k-j}_{s=0}\lt(\begin{array}{c}
  k-j  \\
  s
\end{array}\rt)\lt(\partial^{k-j-s}\overline{\partial}^{k-j}\phi\rt)
\partial^se^{-\vp}\\
&=(-1)^k\sum^{k}_{j=1}A_{k-j}\partial^{k-j}
\lt(\lt(\overline{\partial}^{k-j}\phi\rt)e^{-\vp}\rt)
\end{align*}
Therefore, as the key step of the proof, we have
\begin{align*}
\lt\ll\phi,\overline{\partial}^k\lt(\overline{\partial}_\vp^{k*}\phi\rt)
-\overline{\partial}^{k*}_\vp\lt(\overline{\partial}^k\phi\rt)
\rt\rl_\vp
&=\int_\mb{C}\overline{\phi}\lt(\overline{\partial}^k\lt(\overline{\partial}_\vp^{k*}\phi\rt)
-\overline{\partial}^{k*}_\vp\lt(\overline{\partial}^k\phi\rt)
\rt)e^{-\vp}d\sigma\\
&=(-1)^k\sum^{k}_{j=1}A_{k-j}\int_\mb{C}\overline{\phi}\partial^{k-j}
\lt(\lt(\overline{\partial}^{k-j}\phi\rt)e^{-\vp}\rt)d\sigma\\
&=(-1)^k\sum^{k}_{j=1}A_{k-j}(-1)^{k-j}\int_\mb{C}
\lt(\partial^{k-j}\overline{\phi}\rt)
\lt(\lt(\overline{\partial}^{k-j}\phi\rt)e^{-\vp}\rt)d\sigma\\
&=(-1)^k\sum^{k}_{j=1}A_{k-j}(-1)^{k-j}\int_\mb{C}
\overline{\lt(\overline{\partial}^{k-j}\phi\rt)}
\lt(\overline{\partial}^{k-j}\phi\rt)e^{-\vp}d\sigma\\
&=(-1)^k\sum^{k}_{j=1}A_{k-j}(-1)^{k-j}\lt\ll \overline{\partial}^{k-j}\phi,\overline{\partial}^{k-j}\phi\rt\rl_\vp\\
&=(-1)^k\sum^{k}_{j=1}A_{k-j}(-1)^{k-j}\lt\| \overline{\partial}^{k-j}\phi\rt\|^2_\vp\\
&=\sum^{k}_{j=1}\frac{(k!)^2}{((k-j)!)^2j!}\lt\| \overline{\partial}^{k-j}\phi\rt\|^2_\vp\\
&=\sum^{k-1}_{j=0}\frac{(k!)^2}{(j!)^2(k-j)!}\lt\| \overline{\partial}^{j}\phi\rt\|^2_\vp.
\end{align*}
Then (\ref{28}) is proved.
\ep

\section{Proof of theorems}

First we give the proof of Theorem \ref{thzy}.

\bp
Let $\vp=|z|^2$. By Lemma \ref{lemmaH} and Lemma \ref{lemmak2}, we have for $\forall\phi\in C_0^\infty(\mb{C})$,
\begin{align}
\lt\|\lt(\overline{\partial}^k+a\rt)^*_\vp\phi\rt\|^2_\vp
&\geq\lt\ll\phi,\overline{\partial}^k\lt(\overline{\partial}_\vp^{k*}\phi\rt)
-\overline{\partial}^{k*}_\vp\lt(\overline{\partial}^k\phi\rt)
\rt\rl_\vp\nonumber\\
&=\sum^{k-1}_{j=0}\frac{(k!)^2}{(j!)^2(k-j)!}\lt\| \overline{\partial}^{j}\phi\rt\|^2_\vp\nonumber\\
&\geq k!\lt\|\phi\rt\|^2_\vp.\label{34}
\end{align}
By Cauchy-Schwarz inequality and (\ref{34}), we have for $\forall\phi\in C_0^\infty(\mb{C})$,
\begin{align*}
|\langle f,\phi\rangle_\vp|^2
&\leq\lt\|f\rt\|^2_\vp
\lt\|\phi\rt\|^2_\vp\\
&=\lt(\frac{1}{k!}\lt\|f\rt\|^2_\vp\rt)
\lt(k!\lt\|\phi\rt\|^2_\vp\rt)\\
&\leq\lt(\frac{1}{k!}\lt\|f\rt\|^2_\vp\rt)
\lt\|\lt(\overline{\partial}^k+a\rt)^*_\vp\phi\rt\|^2_\vp.
\end{align*}
Then
by Lemma \ref{lemmaifif}, there exists $u\in L^2(\mb{C},e^{-\vp})$ such that
$$\overline{\partial}^ku+au=f \ \ {\it with} \ \ \|u\|^2_\vp\leq \frac{1}{k!}\lt\|f\rt\|^2_\vp.$$
The proof is complete.
\ep

Second we prove the following theorem.

\bt
There exists a bounded (linear) operator $T_k:L^2(\mb{C},e^{-|z|^2})\rightarrow L^2(\mb{C},e^{-|z|^2})$ such that
$$\lt(\overline{\partial}^k+a\rt)T_k=I \ \ {\it with} \ \ \|T_k\|\leq\frac{1}{\sqrt{k!}},$$
where $\|T_k\|$ is the norm of $T_k$ in $L^2(\mb{C},e^{-|z|^2})$.
\et

\bp
Let $\vp=|z|^2$. For each $f\in L^2(\mb{C},e^{-\vp})$, from Theorem \ref{thzy},
there exists $u\in L^2(\mb{C},e^{-\vp})$ such that
$$\lt(\overline{\partial}^k+a\rt)u=f \ \ {\it with} \ \ \|u\|_\vp\leq\frac{1}{\sqrt{k!}}\lt\|f\rt\|_\vp.$$
Denote this $u$ by $T_k(f)$. Then $T_k(f)$ satisfies $$\lt(\overline{\partial}^k+a\rt)T_k(f)=f \ \ {\it with} \ \ \|T_k(f)\|_\vp\leq\frac{1}{\sqrt{k!}}\lt\|f\rt\|_\vp.$$
Note that $f$ is arbitrary in $L^2(\mb{C},e^{-\vp})$. So $T_k:L^2(\mb{C},e^{-\vp})\rightarrow L^2(\mb{C},e^{-\vp})$ is a bounded (linear) operator such that
$$\lt(\overline{\partial}^k+a\rt)T_k=I \ \ {\it with} \ \ \|T_k\|\leq\frac{1}{\sqrt{k!}}.$$ The proof is complete.
\ep

Lastly we prove Theorem \ref{th1}.

\bp
From Lemma \ref{lemmak}, we have
$$\overline{\partial}\lt(\overline{\partial}_\vp^{*}\phi\rt)-\overline{\partial}^{*}_\vp
\lt(\overline{\partial}\phi\rt)=\phi \overline{\partial}\partial\vp, \ \ \forall\phi\in C_0^\infty(\mb{C}).$$
Then by Lemma \ref{lemmaH}, for $\forall\phi\in C_0^\infty(\mb{C})$,
\begin{align}
\lt\|\lt(\overline{\partial}+a\rt)^*_\vp\phi\rt\|^2_\vp
\geq\lt\langle \phi,\overline{\partial}\lt(\overline{\partial}_\vp^{*}\phi\rt)
-\overline{\partial}^{*}_\vp
\lt(\overline{\partial}\phi\rt)\rt\rangle_\vp
=\lt\langle \phi,\phi \overline{\partial}\partial\vp\rt\rangle_\vp
=\lt\|\phi\sqrt{\overline{\partial}\partial\vp}\rt\|^2_\vp.\label{33}
\end{align}
By Cauchy-Schwarz inequality and (\ref{33}), we have for $\forall\phi\in C_0^\infty(\mb{C})$,
\begin{align*}
|\langle f,\phi\rangle_\vp|^2
&=\lt|\lt\langle \frac{f}{\sqrt{\overline{\partial}\partial\vp}},
\phi\sqrt{\overline{\partial}\partial\vp}\rt\rangle_\vp\rt|^2\\
&\leq\lt\|\frac{f}{\sqrt{\overline{\partial}\partial\vp}}\rt\|^2_\vp
\lt\|\phi\sqrt{\overline{\partial}\partial\vp}\rt\|^2_\vp\\
&\leq\lt\|\frac{f}
{\sqrt{\overline{\partial}\partial\vp}}\rt\|^2_\vp\lt\|
\lt(\overline{\partial}+a\rt)^*_\vp\phi\rt\|^2_\vp.
\end{align*}
Then
by Lemma \ref{lemmaifif}, there exists $u\in L^2(\mb{C},e^{-\vp})$ such that
$$\overline{\partial}u+au=f \ \ {\it with} \ \ \|u\|^2_\vp\leq \lt\|\frac{f}{\sqrt{\overline{\partial}\partial\vp}}\rt\|^2_\vp.$$
The proof is complete.
\ep

\section{Further remarks}

\noindent\textbf{Remark 1.}
Given $\lambda>0$ and $z_0\in\mb{C}$, for the weight $\vp=\lambda|z-z_0|^2$, we obtain the following corollary from Theorem \ref{thzy}. Here we stress that the proof is not simply a straightforward scaling, instead it will scale to a different equation.

\bc\label{coro1}
For each $f\in L^2(\mb{C},e^{-\lambda|z-z_0|^2})$, there exists a weak solution
$u\in L^2(\mb{C},e^{-\lambda|z-z_0|^2})$ solving the equation
$$\overline{\partial}^ku+au=f$$
with the norm estimate $$\int_\mb{C} |u|^2e^{-\lambda|z-z_0|^2}d\sigma\leq
\frac{1}{\lambda^kk!}\int_\mb{C} |f|^2e^{-\lambda|z-z_0|^2}d\sigma.$$
\ec

\bp
From $f\in L^2(\mb{C},e^{-\lambda|z-z_0|^2})$, we have
\begin{equation}\label{1001}
\int_{\mb{C}} |f(z)|^2e^{-\lambda|z-z_0|^2}d\sigma<+\infty.
\end{equation}
Let $z=\frac{\omega}{\sqrt{\lambda}}+z_0$ and $g(\omega)=f(z)= f\lt(\frac{\omega}{\sqrt{\lambda}}+z_0\rt)$. Then by (\ref{1001}), we have
\begin{equation*}
\frac{1}{\lambda}\int_{\mb{C}} |g(\omega)|^2e^{-|\omega|^2}\frac{1}{2i}d\overline{\omega}\wedge d\omega<+\infty,
\end{equation*}
which implies that $g\in L^2(\mb{C},e^{-|\omega|^2})$. For $g$, applying Theorem \ref{thzy} with $a$ replacing by $\frac{a}{(\sqrt{\lambda})^k}$, there exists a weak solution $v\in L^2(\mb{C},e^{-|\omega|^2})$ solving the equation
\begin{equation}\label{1003}
\overline{\partial}^k v(\omega)+\frac{a}{(\sqrt{\lambda})^k}v(\omega)=g(\omega)
\end{equation}
in $\mb{C}$ with the norm estimate
\begin{equation}\label{1004}
\int_{\mb{C}} |v(\omega)|^2e^{-|\omega|^2}\frac{1}{2i}d\overline{\omega}\wedge d\omega\leq
\frac{1}{k!}\int_{\mb{C}} |g(\omega)|^2e^{-|\omega|^2}\frac{1}{2i}d\overline{\omega}\wedge d\omega.
\end{equation}
Note that $\omega=\sqrt{\lambda}(z-z_0)$ and $g(\omega)=f(z)$. Let $u(z)=\frac{1}{(\sqrt{\lambda})^k}v(\omega)=\frac{1}{(\sqrt{\lambda})^k}v\lt(\sqrt{\lambda}(z-z_0)\rt)$. Then (\ref{1003}) and (\ref{1004}) can be rewritten by
\begin{equation}\label{1009}
\overline{\partial}^k u(z) +au(z)
=f(z)
\end{equation}
\begin{equation}\label{1010}
\int_\mb{C} |u(z)|^2e^{-\lambda|z-z_0|^2}d\sigma\leq
\frac{1}{\lambda^kk!}\int_\mb{C} |f(z)|^2e^{-\lambda|z-z_0|^2}d\sigma.
\end{equation}
(\ref{1010}) implies that $u\in L^2(\mb{C},e^{-\lambda|z-z_0|^2})$. Then by (\ref{1009}) and (\ref{1010}), the proof is complete.
\ep

\noindent\textbf{Remark 2.}
From Corollary \ref{coro1}, we can obtain the following corollary, which shows that for any choice of $a$, the differential operator $\overline{\partial}^k+a$ has a bounded right inverse in $L^2(U)$, provided $U$ is a bounded open set.

\bc\label{thzwj1}
Let $U\subset\mb{C}$ be any bounded open set. For each $f\in L^2(U)$, there exists a weak solution $u\in L^2(U)$ solving the equation
$$\overline{\partial}^ku+au=f$$
with the norm estimate $\|u\|_{L^2(U)}\leq c\|f\|_{L^2(U)}$,
where the constant $c$ depends only on the diameter of $U$.
\ec

\bp
Let $z_0\in U$. For given $f\in L^2(U)$, extending $f$ to zero on $\mb{C}\setminus U$, we have
\begin{align*}
\tilde{f}=\left\{
\begin{array}{ccc}
f,       &      & {x\in U}\\
0,     &      & {x\in \mb{C}\setminus U.}
\end{array} \right.
\end{align*}
Then $\tilde{f}\in L^2(\mb{C})\subset L^2(\mb{C},e^{-|z-z_0|^2})$. From Corollary \ref{coro1},
there exists $\tilde{u}\in L^2(\mb{C},e^{-|z-z_0|^2})$ such that
$$\overline{\partial}^k\tilde{u}+a\tilde{u}=\tilde{f} \ \ {\it with} \ \ \int_{\mb{C}}|\tilde{u}|^2e^{-|z-z_0|^2}d\sigma
\leq\frac{1}{k!}\int_{\mb{C}}|\tilde{f}|^2e^{-|z-z_0|^2}d\sigma.$$ Then
\begin{align*}
\int_{\mb{C}}|\tilde{u}|^2e^{-|z-z_0|^2}d\sigma
\leq\frac{1}{k!}\int_{\mb{C}}|\tilde{f}|^2d\sigma
=\frac{1}{k!}\int_U|f|^2d\sigma.
\end{align*}
Note that
\begin{align*}
\int_{\mb{C}}|\tilde{u}|^2e^{-|z-z_0|^2}d\sigma
\geq\int_U|\tilde{u}|^2e^{-|z-z_0|^2}d\sigma
\geq\int_U|\tilde{u}|^2e^{-|U|^2}d\sigma
=e^{-|U|^2}\int_U|\tilde{u}|^2d\sigma,
\end{align*}
where $|U|$ is the diameter of $U$.
Therefore,
$$
e^{-|U|^2}\int_U|\tilde{u}|^2d\sigma
\leq \frac{1}{k!}\int_U|f|^2d\sigma,$$
i.e.,
$$\int_U|\tilde{u}|^2d\sigma
\leq \frac{e^{|U|^2}}{k!}\int_U|f|^2d\sigma.$$
Restricting $\tilde{u}$ on $U$ to get $u$,
then
$$\overline{\partial}^ku+au=f \ \ {\it with} \ \ \int_U|u|^2d\sigma
\leq \frac{e^{|U|^2}}{k!}\int_U|f|^2d\sigma.$$
Note that $u\in L^2(U)$ and let $c=\sqrt{\frac{e^{|U|^2}}{k!}}$. Then the proof is complete.
\ep

\noindent\textbf{Remark 3.}
As a simple consequence of Theorem \ref{thzy}, we can obtain the following
result on the existence of entire weak solutions of the equation $\overline{\partial}^ku+au=f$ for square integrable functions and almost everywhere bounded functions.

\bc\label{thx1}
For each $f\in L^2(\mb{C})$ or $f\in L^\infty(\mb{C})$,
there exists a weak solution $u\in L^2_{loc}(\mb{C})$ solving the equation
$$\overline{\partial}^ku+au=f.$$
In particular, the equation $\overline{\partial}^ku=f$
has a weak solution $u\in L^2_{loc}(\mb{C})$ for $f\in L^2(\mb{C})$ or $f\in L^\infty(\mb{C})$.
\ec

The proof of Corollary \ref{thx1} follows from the observation that $L^2(\mb{C})\subset L^2(\mb{C},e^{-|z|^2})$, $L^\infty(\mb{C})\subset L^2(\mb{C},e^{-|z|^2})$ and $L^2(\mb{C},e^{-|z|^2})\subset L^2_{loc}(\mb{C})$.


\begin{thebibliography}{XX}

\bibitem{1} L. H\"{o}rmander, $L^2$ estimates and existence theorems for the $\overline{\partial}$ operator, Acta Mathematica, 113 (1965), 89-152.

\bibitem{3} L. H\"{o}rmander, An introduction to complex analysis in several variables, D. Van Nostrand Co., Inc., Princeton, N.J.Toronto, Ont.-London, (1966).

\bibitem{4} L. H\"{o}rmander, Notions of convexity, Progress in Mathematics, 127, Birkh\"{a}user Boston, Inc., Boston, MA, (1994).

\bibitem{5} H. Hedenmalm,
On H\"{o}rmander's solution of the $\overline{\partial}$-equation. Math. Z., 281 (1-2), (2015), 349-355.

\bibitem{2} W. P. Johnson, The curious history of Fa\`{a} di Bruno's formula, Amer. Math. Monthly, 109 (3), (2002), 217-234.



\end{thebibliography}
\end{document}